\voffset-10mm
\documentstyle{amsppt}
\magnification1200
\pagewidth{130mm}
\pageheight{200mm}
\hfuzz=2.5pt\rightskip=0pt plus1pt
\binoppenalty=10000\relpenalty=10000\relax
\TagsOnRight
\loadbold
\nologo

\let\wt\widetilde

\let\[\lfloor
\let\]\rfloor

\define\cT{\Cal T}
\define\qbinom#1#2{\thickfracwithdelims[]\thickness0{#1}{#2}_q}
\topmatter
\title
An elementary proof \\
of the irrationality of Tschakaloff series
\endtitle
\author
Wadim Zudilin\footnotemark"$^\ddag$"\ {\rm(Moscow)}
\endauthor
\date
E-print \tt math.NT/0506086
\enddate
\address
\hbox to70mm{\vbox{\hsize=70mm%
\leftline{Department of Mechanics and Mathematics}
\leftline{Moscow Lomonosov State University}
\leftline{Vorobiovy Gory, GSP-2}
\leftline{119992 Moscow, RUSSIA}
\leftline{{\it URL\/}: \tt http://wain.mi.ras.ru/}
}}
\endaddress
\email
{\tt wadim\@ips.ras.ru}
\endemail
\abstract
We present a new proof of the irrationality
of values of the series
$\cT_q(z)=\allowmathbreak\sum_{n=0}^\infty z^nq^{-n(n-1)/2}$
in both qualitative and quantitative forms.
The proof is based on a hypergeometric construction
of rational approximations to~$\cT_q(z)$.
\endabstract
\subjclass
11J72, 11J82 (Primary), 33D15 (Secondary)
\endsubjclass
\dedicatory
To A.\,B.~Shidlovski\v\i\ on the occasion of his 90th birthday
\enddedicatory
\endtopmatter
\leftheadtext{W.~Zudilin}
\rightheadtext{Irrationality of Tschakaloff series}
\footnotetext"$^\ddag$"{The work is partially supported
by grant no.~03-01-00359 of the Russian Foundation for Basic Research.}
\document

\subhead
1. Introduction
\endsubhead
In 1919, L.~Tschakaloff introduced the series~\cite{10}
$$
\cT_q(z)=\sum_{n=0}^\infty z^nq^{-n(n-1)/2},
\tag1
$$
convergent in the whole complex $z$-plane whenever $|q|>1$,
and proved the irrationality and linear independence
of its values at rational non-zero points $z$ and $q$
(under certain assumptions on~$q$). His method generalized
that by O.~Sz\'asz~\cite{9} for a special case of~\thetag{1},
namely, the function
$\Theta_q(z)=\sum_{n=0}^\infty z^nq^{-n^2}=\cT_{q^2}(z/q)$;
at about the same time
F.~Bernstein and O.~Sz\'asz~\cite{1} used
a continued fraction for $\Theta_q(z)$ due to Eisenstein
to provide another irrationality proof for its values
at certain rational $q$ and~$z$. These seem to be the very first
results on the arithmetic nature of values
of $q$-series.

The aim of this note is to give an elementary proof
of Tschakaloff's theorem~\cite{10} and also its quantitative
form given by P.~Bundschuh in~\cite{3, Satz~2}.

\proclaim{Theorem}
Let $q=q_1/q_2$ and $z$ be non-zero rational numbers,
where $|q|>1$ and $q_1,q_2\in\Bbb Z$. Suppose that
the non-negative number
$$
\gamma=\frac{\log|q_2|}{\log|q_1|}
$$
satisfies $\gamma<\gamma_0=(3-\sqrt5)/2$. Then
the value $\cT_q(z)$~is irrational.
Moreover, for any $\varepsilon>0$ there exists
a positive constant $b_0(\varepsilon)$ such that
$$
\biggl|\cT_q(z)-\frac ab\biggr|
>|b|^{-1-\frac{\scriptstyle\sqrt5-1}{\scriptstyle2(\gamma_0-\gamma)}-\varepsilon}
\tag2
$$
for all integers $a$ and $b$ with
$|b|\ge b_0(\varepsilon)$.
\endproclaim

The rational approximations to the Tschakaloff function~\thetag{1}
that we construct in the next section are actually
the same as those in~\cite{10}
and~\cite{3}. Our contribution here is to provide an
elementary explanation of why these approximations are good
enough to obtain the irrationality of~$\cT_q(z)$.
Our proof is inspired by the ideas of L.~Gutnik
and Yu.~Nesterenko \cite{7, Section~1}
in their proof that $\zeta(3)\notin\Bbb Q$.
This is the famous theorem due
to R.~Ap\'ery; elementary proofs and interrelations
with irrationality results for other mathematical
constants may be found in~\cite{4} and~\cite{8}.

\subhead
2. Proof
\endsubhead
For the first paragraph, we shall think of~$q$
as a {\it variable}.
Let $n$~be a positive integer and define the polynomial
$$
R(T;q)=R_n(T;q)=(1-qT)(1-q^2T)\dotsb(1-q^nT).
$$
Multiplication gives
$$
R(T;q)
=\sum_{k=0}^nC_k(q)T^k,
\tag3
$$
where, for $k=0,1,\dots,n$,
$$
C_k(q)=C_{k,n}(q)\in\Bbb Z[q]
\tag4
$$
is a polynomial in~$q$ with
$$
\operatorname{degree}C_k(q)\le\frac{n(n+1)}2.
\tag5
$$

Conditions~\thetag{4} and~\thetag{5} imply that, if $q=q_1/q_2$, then
$$
q_2^{n(n+1)/2}C_k\biggl(\frac{q_1}{q_2}\biggr)
\in\Bbb Z
\tag6
$$
for $k=0,1,\dots,n$ and arbitrary non-zero integers $q_1$ and $q_2$.

Let $m=\lfloor\beta n\rfloor$ (here $\lfloor\,\cdot\,\rfloor$~denotes
the integer part of a number), where $\beta=(\sqrt5-1)/2$
is the positive root of the polynomial $x^2+x-1$,
and introduce the series
$$
I_n=I_n(z;q)
=\sum_{t=1}^\infty R_n(q^{-t};q)z^{t+m}q^{-(t+m)(t+m-1)/2},
\tag7
$$
which converges if $|q|>1$.
Using \thetag{3} we obtain
$$
\align
I_n
&=\sum_{t=1}^\infty z^{t+m}\sum_{k=0}^nC_k(q)q^{-kt-(t+m)(t+m-1)/2}
\\
&=\sum_{k=0}^nz^{-k}C_k(q)q^{k(k-1)/2+km}
\sum_{t=1}^\infty z^{k+t+m}q^{-(k+t+m)(k+t+m-1)/2}
\allowdisplaybreak
&=\sum_{k=0}^nz^{-k}C_k(q)q^{k(k-1)/2+km}
\sum_{l=k+m+1}^\infty z^lq^{-l(l-1)/2}
\allowdisplaybreak
&=\sum_{k=0}^nz^{-k}C_k(q)q^{k(k-1)/2+km}
\biggl(\sum_{l=0}^\infty z^lq^{-l(l-1)/2}
-\sum_{l=0}^{k+m}z^lq^{-l(l-1)/2}\biggr)
\\
&=\sum_{k=0}^nz^{-k}C_k(q)q^{k(k-1)/2+km}\cdot\cT_q(z)
-\sum_{k=0}^nC_k(q)\sum_{l=0}^{k+m}z^{-(k-l)}q^{k(k-1)/2+km-l(l-1)/2}.
\tag8
\endalign
$$

If $q=q_1/q_2$ and $z=z_1/z_2$, where
$q_1,q_2,z_1,z_2\in\Bbb Z\setminus\{0\}$, then
from~\thetag{6} and~\thetag{8} we see that
the quantity $\wt I_n=\wt I_n(z;q)$ defined by
$$
\wt I_n
=z_1^nz_2^mq_1^{m(m-1)/2}q_2^{n(n+1)/2+n(n-1)/2+nm}I_n
\tag9
$$
is of the form
$$
\wt I_n=B_n\cdot\cT_q(z)-A_n,
\tag10
$$
where $A_n$ and $B_n$ are integers, determined by~\thetag{8}
and~\thetag{9}.
In addition, since equality in~\thetag{5}
is achieved only when $k=n$, we see that the coefficient
of~$\cT_q(z)$ in~\thetag{8} has the following asymptotics
as $n\to\infty$ (where $f(n)\sim g(n)$ means that
$f(n)/g(n)\to1$):
$$
\align
\biggl|\sum_{k=0}^nz^{-k}C_k(q)q^{k(k-1)/2+km}\biggr|
&\sim|z|^{-n}|C_n(q)|\,|q|^{n(n-1)/2+nm}
\\
&=|z|^{-n}|q|^{n(n+m)}.
\tag11
\endalign
$$

In order to evaluate the asymptotic behavior of
the sum of the series \thetag{7},
notice that $R_n(q^{-t};q)=0$ for $t=1,2,\dots,n$.
Therefore (using $f(n)=O(g(n))$ as $n\to\infty$
to mean that $|f(n)|\le C|g(n)|$
for some constant $C>0$ and all $n$ sufficiently large),
$$
\align
I_n
&=\sum_{t=n+1}^\infty R_n(q^{-t};q)z^{t+m}q^{-(t+m)(t+m-1)/2}
\\
&=R_n(q^{-(n+1)};q)z^{n+m+1}q^{-(n+m)(n+m+1)/2}
+O(q^{-(n+m+1)(n+m+2)/2})
\\
&=z^{n+m+1}q^{-(n+m)(n+m+1)/2}(1-q^{-1})(1-q^{-2})\dotsb(1-q^{-n})
\\ &\qquad
+O(q^{-(n+m+1)(n+m+2)/2})
\\
&\sim z^{n+m+1}q^{-(n+m)(n+m+1)/2}
\qquad\text{as}\quad n\to\infty.
\tag12
\endalign
$$
In particular, $I_n\ne0$ for all $n$ sufficiently large.

Finally, since $|q_1/q_2|=|q|>1$ implies that $|q_1|>1$,
we may define $\gamma$ by the relation
$\log|q_2|=\allowmathbreak\gamma\log|q_1|$, so that
$\gamma\ge0$. Assume that $\gamma<\gamma_0=(3-\sqrt5)/2$.
Then, from \thetag{9}, \thetag{11}, \thetag{12},
and the relation $m=\lfloor\beta n\rfloor$, for the
quantities $B_n$ and $\wt I_n$ in \thetag{10}
we have
$$
\lim_{n\to\infty}\frac{\log|B_n|}{n^2\log|q_1|}
=(1-\gamma)(1+\beta)+\gamma(1+\beta)+\frac{\beta^2}2
=\frac{\sqrt5(\sqrt5+1)}4
\tag13
$$
and
$$
\lim_{n\to\infty}\frac{\log|\wt I_n|}{n^2\log|q_1|}
=-(1-\gamma)\frac{(1+\beta)^2}2+\gamma(1+\beta)+\frac{\beta^2}2
=-\frac{\sqrt5(\sqrt5+1)(\gamma_0-\gamma)}{2(\sqrt5-1)}
<0.
\tag14
$$

Now let us show that $\cT_q(z)$ cannot be rational.
Suppose, on contrary, that $\cT_q(z)=a/b$
for some integers $a$ and $b\ne0$. Then from~\thetag{10}
$$
b\wt I_n=B_na-A_nb\in\Bbb Z
\qquad (n=1,2,\dots).
$$
Recalling that \thetag{12} yields $I_n\ne0$ for $n$~large, we
conclude that $|b\wt I_n|\ge1$. But, by~\thetag{14},
we have $|b\wt I_n|\to0$ as $n\to\infty$.
The contradiction implies that $\cT_q(z)\notin\Bbb Q$.

We leave to the reader the derivation of estimate~\thetag{2}
from~\thetag{10}, \thetag{13}, and~\thetag{14} by letting
$a_n=A_n$ and $b_n=B_n$ in
the following standard lemma (compare
\cite{2, Section~11.3, Exercise~3}).

\proclaim{Lemma}
Let $\alpha$ be an irrational real number. Suppose that we have
a sequence of rational approximations $a_n/b_n$ to~$\alpha$
\rom(where $a_n,b_n\in\Bbb Z$ for $n=1,2,\dots$\rom)
such that the sequence $|b_n|$ tends to infinity with~$n$,
$$
\lim_{n\to\infty}\frac{\log|b_{n+1}|}{\log|b_n|}=1,
$$
and with some constant~$c>0$
$$
\biggl|\alpha-\frac{a_n}{b_n}\biggr|
<\frac1{|b_n|^{1+c}}
$$
for all $n$ sufficiently large.
Then for any $\varepsilon>0$ there exists
a positive constant $b_0(\varepsilon)$ such that
$$
\biggl|\alpha-\frac ab\biggr|
>\frac1{b^{1+1/c+\varepsilon}}
$$
for all integers $a$ and $b$ with
$b\ge b_0(\varepsilon)$.
\endproclaim

\subhead
3. Related results
\endsubhead
Although we are able to prove the irrationality of $\cT_q(z)$
only under the hypothesis $\gamma<\gamma_0=0.381966\dots$, it is expected
that this hypothesis can be dropped, i.e., that $\cT_q(z)$~is
irrational for all $z\in\Bbb Q\setminus\{0\}$ and $q\in\Bbb Q$ with $|q|>1$.
This remains an open problem. The earlier method in~\cite{9} requires
the condition $\gamma<1/3$ (which is worse, since $1/3<\gamma_0$) corresponding
to the simpler choice $\beta=0$ in our notation.
The choice $\beta=(\sqrt5-1)/2$ ensures the optimal value
of~$\gamma_0$ in terms of the construction presented here.

The Tschakaloff function~\thetag{1} might be viewed as
``half'' of the theta series
\linebreak
$\sum_{n\in\Bbb Z}z^{n-1/2}q^{-(n-1/2)^2}$.
This viewpoint and Nesterenko's theorem~\cite{6} on the transcendence
of certain theta series imply the transcendence
of~$\cT_q(z)$ for $q$~algebraic, $|q|>1$,
and $z=q^k$ with some $k\in\Bbb Z$,
solving this case of the open problem.
On the other hand, when $z$ and $q$ are multiplicatively
independent, no transcendence results are
known. This is part of a general problem posed by K.~Mahler
in~\cite{5} for analytic functions which satisfy functional
equations (such as $\cT_q(z)=1+z\cT_q(z/q)$ for the function~\thetag{1}),
but to which his method from~\cite{5} cannot be applied.

The constants $\beta=(\sqrt5-1)/2$ and $\gamma_0=1-\beta$, involved
in the proof of the Theorem, are related to the {\it golden
mean\/} (or {\it golden section\/}),
the positive root of the polynomial $x^2-x-1$.
It is quite curious that the golden mean and its
generalizations (the so-called {\it metallic means\/}) also occur
in other irrationality proofs related to Ap\'ery's theorem~\cite{4}.

Finally, we mention that a special case of the $q$-binomial theorem
implies the following explicit formula for the
polynomial~\thetag{4}:
$$
C_k(q)=(-1)^k\qbinom nkq^{k(k+1)/2}
$$
involving the $q$-binomial coefficients
$$
\qbinom nk=\frac{[n]_q!}{[k]_q!\,[n-k]_q!}\in\Bbb Z[q],
$$
where $[0]_q!=1$ and, for $k=1,2,\dots$,
$$
[k]_q!=\frac{(q-1)(q^2-1)(q^3-1)\dotsb(q^k-1)}{(q-1)^k}.
$$

\subhead
Acknowledgments
\endsubhead
It is a pleasure for me to thank Jonathan Sondow, who
conceptually influenced the note by several very
useful suggestions.



\Refs

\ref\no1
\by F.~Bernstein and O.~Sz\'asz
\paper \"Uber Irrationalit\"at unendlicher
Kettenbr\"uche mit einer Anwendung auf die Reihe
$\dots$
\jour Math. Ann.
\vol76
\yr1915
\pages295--300
\endref

\ref\no2
\by J.\,M.~Borwein and P.\,B.~Borwein
\book Pi and the AGM: A study in analytic number theory
and computational complexity
\bookinfo Canad. Math. Soc. Ser. Monogr. Adv. Texts
\publaddr New York
\publ Wiley
\yr1987
\endref

\ref\no3
\by P.~Bundschuh
\paper Versch\"arfung eines arithmetischen Satzes von Tschakaloff
\jour Portugal. Math.
\vol33
\issue1
\yr1974
\pages1--17
\endref

\ref\no4
\by D.~Huylebrouck
\paper Similarities in irrationality proofs for $\pi$, $\ln2$,
$\zeta(2)$ and $\zeta(3)$
\jour Amer. Math. Monthly
\vol108
\issue3 (March)
\yr2001
\pages222--231
\endref

\ref\no5
\by K.~Mahler
\paper Remarks on a paper by W.~Schwarz
\jour J. Number Theory
\vol1
\yr1969
\pages512--521
\endref

\ref\no6
\by Yu.~Nesterenko
\paper Modular functions and transcendence problems
\jour C.~R. Acad. Sci. Paris S\'er.~I
\vol322
\issue10
\yr1996
\pages909--914
\endref

\ref\no7
\by Yu.\,V.~Nesterenko
\paper A few remarks on~$\zeta(3)$
\jour Math. Notes
\vol59
\issue6
\yr1996
\pages625--636
\endref

\ref\no8
\by A.~van der Poorten
\paper A proof that Euler missed...
Ap\'ery's proof of the irrationality of~$\zeta(3)$
\paperinfo An informal report
\jour Math. Intelligencer
\vol1
\issue4
\yr1978/79
\pages195--203
\endref

\ref\no9
\by O.~Sz\'asz
\paper \"Uber Irrationalit\"at gewisser unendlicher Reihen
\jour Math. Ann.
\vol76
\yr1915
\pages485--487
\endref

\ref\no10
\by L.~Tschakaloff
\paper Arithmetische Eigenschaften der unendlichen Reihe
$\dots$
\jour Math. Ann.
\vol80
\yr1919
\pages62--74
\endref

\endRefs
\enddocument